\documentclass[reqno]{article}

\usepackage{xcolor}
\usepackage{amssymb}
\usepackage{amsfonts}
\usepackage{amsmath}
\numberwithin{equation}{section}
\usepackage{enumerate}
\newtheorem{theorem}{Theorem}[section]
\newtheorem{proposition}[theorem]{Proposition}
\newtheorem{lemma}[theorem]{Lemma}
\newtheorem{corollary}[theorem]{Corollary}
\newtheorem{definition}[theorem]{Definition}
\newtheorem{example}[theorem]{Example}

\newcommand{\dom}{{\rm dom}}

\title {Novel approach to root functions of matrix polynomials with applications in differential equations and meromorphic matrix functions} 
\author{Author: Muhamed Borogovac}

\begin{document}
\maketitle 
\textbf{Abstract} 
In the first part of the paper, we address an invertible matrix polynomial $L(z)$ and its inverse $\hat{L}(z) := -L(z)^{-1}$. We present a method for obtaining a canonical set of root functions and Jordan chains of $L(z)$ through elementary transformations of the matrix $L(z)$ alone. This method provides a new and simple approach to deriving a general solution of the system of ordinary linear differential equations $L\left(\frac{d}{dt}\right)u=0$ using only elementary transformations of the corresponding matrix polynomial $L(z)$. 

In the second part of the paper, given a matrix generalized Nevanlinna function $Q\in N_{\kappa }^{n \times n}$ and a canonical set of root functions of $\hat{Q}(z) := -Q(z)^{-1}$, we provide an algorithm to determine a specific Pontryagin space $(\mathcal{K}, [.,.])$, a specific self-adjoint operator $A:\mathcal{K}\rightarrow \mathcal{K}$ and an operator $\Gamma: \mathbb{C}^{n}\rightarrow \mathcal{K}$ that represent the function $Q$ in a Krein-Langer type representation. We demonstrate the main results through examples of linear systems of ODEs.



\section{Preliminaries and introduction}\label{s2}
\textbf{1.1} Let $\mathbb{N}$, $\mathbb{R}$, $\mathbb{C}$ denote sets of positive integers, real numbers, and complex 
numbers, respectively. We assume that the reader is familiar with basic topics of Functional Analysis, particularly Hilbert spaces.  Let $(\cdot, \cdot)$ denote the definite scalar product in a Hilbert space $\mathcal{H}$, and let $\mathcal{L}(\mathcal{H})$ denote the Banach space of bounded linear operators on $\mathcal{H}$. Let $(\mathcal{K}, [\cdot, \cdot])$ denote a Krein space, which is a complex vector space with a scalar product, i.e., a Hermitian sesquilinear form $[\cdot, \cdot]$, defined such that the following decomposition of $\mathcal{K}$ exists:
\[
\mathcal{K}=\mathcal{K}_{+}[+]\mathcal{K}_{-},
\]
where $(\mathcal{K}{+},[\cdot, \cdot])$ and $(\mathcal{K}_{-},-[\cdot, \cdot])$ are Hilbert spaces that are mutually orthogonal with respect to the form $[\cdot, \cdot]$. For every Krein space $(\mathcal{K}, [\cdot, \cdot])$, one can associate a Hilbert space $(\mathcal{K}, (\cdot, \cdot))$ by
\[
\left( x,y \right)=\left[ Jx,y \right],\, \forall x,y\in\mathcal{K},
\]
where $J$ is an operator that satisfies $J^{-1}=J^{\ast}=J$. This operator is called the \textit{fundamental symmetry}. If the subspace $\mathcal{K}_{-}$ is $\kappa-$dimensional, where $\kappa \in \mathbb{N}$, then we call it a \textit{Pontryagin space} of index $\kappa$. For properties of Pontryagin spaces, see \cite{IKL}.
\\

\textbf{1.2.} Given $n\times n$ constant complex matrices 
$A_{j},\, \, j=0,\, 1,\, \mathellipsis ,\, l$, where $l \in \mathbb{N}$, and the corresponding linear 
system of differential equations:
\begin{equation}
\label{eq16}
A_{l}\frac{d^{l}u}{{dt}^{l}}+\mathellipsis +A_{1}\frac{du}{dt}+A_{0}u=f\left( t \right),
\end{equation}
where $u(t):=\left(
\begin{array}{*{20}c}
u_{1}(t)\\
\mathellipsis \\ 
u_{n}(t) \\ 
\end{array} \right)$ is an n-dimensional unknown vector function and $f\left( t \right)$ is a piecewise continuous $n-$dimensional vector function defined on the real axis, i.e., on $\mathbb{R}$. This is called a \textit{linear system of ordinary differential equations (ODE) with constant coefficients}. If we introduce operator $L\left( \frac{d}{dt} \right)$ by:
\[
L\left( \frac{d}{dt} \right)u:=A_{l}\frac{d^{l}u}{{dt}^{l}}+\mathellipsis 
+A_{1}\frac{du}{dt}+A_{0}u,
\]
then 
\begin{equation}
\label{eq18}
L\left( \frac{d}{dt} \right)u=0
\end{equation}
is the homogeneous system associated with (\ref{eq16}). One can seek a solution in the 
form
\begin{equation}
\label{eq110}
u\left( z \right)=\varphi e^{zt},
\end{equation}
where $\varphi :=\left( {\begin{array}{*{20}c}
\varphi_{1}\\
\mathellipsis \\
\varphi_{n}\\
\end{array} } \right)\ne 0,\, \varphi_{i}\in \mathbb{C},\, i=1,2,\, \mathellipsis ,\, n$. After substituting (\ref{eq110}) into (\ref{eq18}), we obtain:
\begin{equation}
\label{eq112}
L\left( z \right)\varphi e^{zt}=0,
\end{equation}
where 
\[
L\left( z \right):=A_{l}z^{l}+\mathellipsis 
+A_{1}z+A_{0} .
\]
The matrix polynomial $L(z)$ is called \textit{monic} if $A_{l}=I$. If an n-dimensional vector $\varphi \ne 0$ is a solution of (\ref{eq112}) for some $z=\alpha$, i.e., a solution of the algebraic system:
\begin{equation}
\label{eq114}
L\left( \alpha \right)\varphi=0,
\end{equation}
then $z=\alpha$ is a solutions of the equation:
\begin{equation}
\label{eq116}
\det {L\left( z \right)}=0.
\end{equation}
The polynomial $\chi (z) :=\det {L\left( z \right)}$ is called the \textit{characteristic polynomial} of the matrix polynomial $L(z)$. The zeros of the characteristic polynomial, i.e., solutions of the equation (\ref{eq116}) are called \textit{eigenvalues} of $L\left( z \right)$. Every solution $\varphi \ne 0$ of (\ref{eq114}) is called an \textit{eigenvector} corresponding to the eigenvalue $\alpha $. The set of all eigenvalues of $L\left( z \right)$ is called the finite spectrum of $L\left( z \right)$ and is denoted by $\sigma_{p}\left( L \right)$ or simply $\sigma \left( L \right)$. 

The following proposition, which gives a general solution of the homogeneous system (\ref{eq18}), is an equivalent version of \cite[Proposition 1.9]{GLR}. We state it here for the convenience of the reader. 

\begin{proposition}\label{proposition14} The vector function 
\begin{equation}
\label{eq118}
u\left( t \right)=\left( \frac{t^{k-1}}{(k-1)!}\varphi 
_{0}+\frac{t^{k-2}}{(k-2)!}\varphi_{1}+\mathellipsis +\varphi_{k-1} 
\right)e^{\alpha t},
\end{equation}
where $k\in \mathbb{N}$ and  $\varphi_{j}\in \mathbb{C}^{n}$ for $j=0,\, 1,\, \mathellipsis ,\, k-1$, is a solution of equation (\ref{eq18}) if and only if the following equalities hold: 
\begin{equation}
\label{eq120}
\sum\limits_{p=0}^i {\frac{1}{p!}L^{\left( p \right)}\left( \alpha 
\right)\varphi_{i-p}} =0, \, i=0,\, \mathellipsis ,\, k-1,
\end{equation}
where $L^{\left( p \right)}\left( \alpha \right)$ is p-th derivative 
of $L\left( z \right)$ at the eigenvalue $\alpha $ of $L\left( z \right)$.
\end{proposition}
The sequence of n-dimensional vectors $\varphi_{0},\, \varphi_{1},\, \mathellipsis ,\, \varphi_{k-1}$, $\varphi_{0}\ne 0$, for which identities (\ref{eq120}) hold, is called a \textit{Jordan chain of order or length} $k$ \textit{for} $L\left( z \right)$ \textit{at the eigenvalue} $\alpha \in \mathbb{C}$. The Jordan chain is said to be \textit{maximal, of order} $k$, if the additional equation, for $i=k$, in the system (\ref{eq120}) does not have a solution. The \textit{canonical system} of Jordan chains for the eigenvalue $\alpha $ of $L\left( z \right)$ consists of all linearly independent eigenvectors at $\alpha $ and Jordan chains of maximal length, one for each of the eigen-vectors. The linear span of these Jordan vectors at $\alpha $ is called the \textit{algebraic eigenspace} of $L(z)$ at $\alpha $. The \textit{canonical system of} $L\left( z \right)$ is the union of canonical systems of $L(z)$ over all eigenvaues of $L(z)$. 

It is easy to see that all above definitions apply to the Jordan chain $f_{0},\, f_{1},\, \mathellipsis ,\, f_{k-1}$ with $ f_{0}\ne 0$ at the eigenvalue $\alpha \in \mathbb{C}$ of a constant square matrix $A$. Namely, in that case $L(z)=A-zI$.

In \cite[Section 1.6]{GLR}, authors provide one algorithm for finding a canonical Jordan chain of any matrix polynomial $L\left( z \right)$ at $\alpha $. In this paper, we describe how to find a canonical system of Jordan vectors at $\alpha $  using only elementary transformations of the matrix polynomial $L\left( z \right)$. In Example \ref{example214}, we demonstrate how to apply this method to find Jordan vectors of a given invertible matrix polynomial and solve the corresponding system (\ref{eq18}) of differential equations. 
\\

\textbf{1.3.} Recall that an operator valued function $Q:D\left( Q \right)\subset \mathbb{C}\to \mathcal{L}(\mathcal{H})$ belongs to the \textit{generalized Nevanlinna class} $N_{\kappa }\left( \mathcal{H} \right)$ if it is meromorphic on $\mathbb{C} \backslash \mathbb{R}$, satisfies ${Q\left( z \right)}^{\ast }=Q\left(\bar{z} \right)$ for all points $z$ of holomorphy of $Q$, and the kernel $N_{Q}\left( z,w \right):=\frac{Q\left( z \right)-{Q\left( w \right)}^{\ast 
}}{z-\bar{w}}$ has $\kappa $ negative squares. Then $\kappa \in \mathbb{N}\cup\lbrace 0 \rbrace$ is called 
a negative index of $Q$. 

Representations of generalized Nevanlinna scalar and matrix functions, denoted by $Q \in {N_{\kappa}}^{n \times n}$, in terms of symmetric (Hermitian) operators in Pontryagin spaces were developed by M. G. Krein and H. Langer for functions $Q$ that satisfy certain properties at infinity  (see \cite{KL1, KL2}). Those operator  representations are frequently called \textit{Krein-Langer} representations (see \cite{HSW}). These representations were later generalized to all functions $Q \in {N_{\kappa}}^{n \times n}$ by means of linear relations (see e.g. \cite{DLS}). For properties of linear relations, i.e., multi-valued linear operators, one can refer to \cite{A,Lu1}. The most general representation, when $Q \in N_{\kappa }\left( \mathcal{H} \right)$, is provided below for the reader's convenience. 

\begin{theorem}\label{theorem12} \cite[Proposition 2.3 ]{Lu3} A function $Q:\dom \, Q \to \mathcal{L}(\mathcal{H})$ is a generalized Nevanlinna function if and only if there exist a Pontryagin space $(\mathcal{K}, [.,.])$, a self-adjoint relation A in  $\mathcal{K}$, a point $z_{0} \in \rho(A)\cap \mathbb{C}^{+}$ and a bounded linear map  $\Gamma :\mathcal{H}\to \mathcal{K}$ such that $Q$ can be written as
\begin{equation}
\label{eq12}
Q\left( z \right)={Q(z_{0})}^{\ast }+(z-\bar{z_{0}})\Gamma^{+}\left(I+\left( z-z_{0} \right)\left( A-z \right)^{-1} \right)\Gamma ,\, z\in \dom\,\left( Q \right),
\end{equation}
Moreover, this realization can be chosen minimal, that is,
\[
\mathcal{K}=c.l.s.\left\{ \Gamma_{z}h:z\in \dom\,Q,h\in \mathcal{H} \right\}
\]
where
\begin{equation}
\label{eq14}
\Gamma_{z}=\left( I+\left( z-z_{0} \right)\left( A-z \right)^{-1} \right)\Gamma .
\end{equation}
In this case the realization (\ref{eq12}) is unique up to unitary equivalence and $Q \in N_{\kappa}\left( \mathcal{H}\right) $ if and only if index of the Pontryagin space equals $\kappa$.
\end{theorem}

When $L\left( z\right)$ is a Hermitian matrix polynomial with $det\, L\left( z \right)\not\equiv 0$, then $L\left( z\right) $ is a generalized Nevanlinna matrix function, i.e., $L \in {N_{\kappa}}^{n \times n}$.  In that case, the inverse function $\hat{L}\left( z\right):=-{L\left( z\right)}^{-1}$ also satisfies $\hat{L} \in {N_{\kappa}}^{n \times n}$. Therefore, Theorem \ref{theorem12}, and a special case of it, Lemma \ref{lemma32} below, applies to $\hat{L}(z)$.

It is well known that the generalized poles (including poles) of a generalized Nevanlinna function $Q\in N_{\kappa}(\mathcal{H})$, represented by the minimal operator representation (\ref{eq12}) coincide with the eigenvalues of the representing operator $A$. Moreover, in recent references, generalized poles and their orders are defined in terms of the eigenvalues of $A$ and the orders of the corresponding canonical Jordan chains, i.e., in terms of root subspaces of $A$ (see \cite[Definition 2.4]{BLu}, \cite[Definition 4.2]{Lu3}).

Therefore, the generalized poles, including their orders, of a function $Q$ can be investigated by examining the eigenvalues and Jordan chains of the representing operator $A$, and vice-verse. This makes operator representations a very important tool in Functional analysis. However, it is usually very difficult to find the representing operators $A$ and $\Gamma$ in the Krein-Langer operator representation (\ref{eq12}) of a given function $Q$. The constructions used in the cited papers are abstract and not applicable in concrete situations. 

The following is a brief review of the results in Section \ref{s6}. In Proposition \ref{proposition34}, we provide a special case of the Krein-Langer representation of generalized Nevanlinna functions. This representation is applicable to the inverse functions $\hat{L}(z)=-L(z)^{-1}$. We use this result in Corollary \ref{corollary310} to give an operator representation of the function $\hat{L}(z)$ that vanishes at infinity. 

In Theorem \ref{theorem38}, we prove a method that allows us to construct a specific Krein space $(\mathcal{K}, [.,.])$, a specific self-adjoint operator $A$ in $\mathcal{K}$, and an operator $\Gamma: \mathbb{C}^{n}\rightarrow \mathcal{K}$ to represent the given function $Q \in N^{n \times n}$. In Example \ref{example312}, we apply the results of Theorem \ref{theorem38} to a given system of ODEs with a matrix polynomial $L(z)$ to obtain the Krein-Langer representation corresponding to $\hat{L}(z)$. In this example, we first use Proposition \ref{proposition26} to find the canonical set of root functions of $L(z)$. Then, we demonstrate how to use a canonical set of root functions to find the operator representation of $\hat{L}(z)$ by the method proved in Theorem \ref{theorem38}, when $\sigma(L)\subset \mathbb{R}$.



\section{Finding a root function of $L\left( z \right)$ by means of elementary matrix transformations}\label{s4}
\textbf{2.1.} Let $L\left( z \right)$ be an $n\times n$ invertible matrix polynomial. This implies $det\, L\left( z \right)\not\equiv 0$. It is well known that, by means of elementary matrix transformations of $ L\left( z \right)$, we can obtain the following matrix $D\left( z \right)$:
\begin{equation}
\label{eq22}
D\left( z \right)=S\left( z \right)L\left( z \right)T\left( z \right)
\end{equation}
with
\begin{equation}
\label{eq24}
D\left( z \right):=\left[ {\begin{array}{*{20}c}
d_{1}\left( z \right) & \cdots & 0\\
\vdots & \ddots & \vdots \\
0 & \cdots & d_{n}\left( z \right)\\
\end{array} } \right],
\end{equation}
where $d_{i}(z)$ are polynomials, and $S\left( z \right)$ and $T\left( z \right)$ are $n\times n$ matrix 
polynomials, with $det\, S\left( z \right)$ and 
$det\, T(z)$ being constant non-zero numbers. 
The matrix $D\left( z \right)$ in (\ref{eq24}) is referred to as the \textit{diagonal form of the matrix polynomial} $L\left( z \right)$, and $S\left( z \right)$ and $T\left( z \right)$ as the \textit{matrices of elementary transformations}. Note that the product of polynomials $d_{i}(z)$ is uniquely determined, while their order on the diagonal is not; it depends on the sequence of the elementary transformations, i.e., on the matrices $S(z)$ and $T(z)$.

For a proof of the above statement see, for example, \cite[Theorem S1.1]{GLR}. Note that the algorithm in that proof yields the matrix $D(z)$ with additional properties: the polynomials $d_{i}(z)$ are monic, and each $d_{i}(z)$ is divisible by $d_{i-1}(z)$. In this case, the matrix $D(z)$ is called the \textit{Smith form} of the matrix polynomial $L(z)$. Because we do not need the Smith form, in many concrete situations we may avoid the routine given in the proof of \cite[Theorem S1.1]{GLR} and find a more efficient sequence of elementary transformations to diagonalize $L(z)$.

If for a fixed $i \in \lbrace 1, ..., n\rbrace $ we find all $\beta_{i}$ solutions $\alpha_{i,j}, j=1, ..., \beta_{i}$, of the equation $d_{i}\left( z \right)=0$ with their multiplicities, or orders, $k_{i,j}$, then it holds
\begin{equation}
\label{eq26}
d_{i}\left( z \right)=\left( z-\alpha_{i,1} \right)^{k 
_{i1}}\mathellipsis \left( z-\alpha_{i,\beta_{i}} \right)^{k 
_{{i\beta}_{i}}},\, i=1,\mathellipsis ,n.
\end{equation} 
It is possible to have equal factors $\left( z-\alpha \right)$ in the decompositions (\ref{eq26}) of different polynomials $d_{i}\left( z \right)$, but we still consider those factors as different because they are in different places in the matrix polynomial $D\left( z \right)$ and correspond to different eigenvectors. They also may have different orders $k_{i}$ in different polynomials $d_{i}$.

Let us recall the easiest way to obtain matrices $S\left( z \right)$ and $T\left( z \right)$. We start with the block matrix $\left( {\begin{array}{*{20}c}
I_{11} & L(z)\\
0 & I_{22}\\
\end{array} } \right)$, where $I_{ii},\, i=1,2$, denote $n \times n$ identity matrices. Then we perform all necessary elementary transformations of rows and columns of the block matrix so that $L(z)$ transforms to a diagonal form $D(z)$. In other words, we perform the same transformations on the rows of the matrices $L(z)$ and $I_{11}$ simultaneously, and the same transformations on the columns of the matrices $L(z)$ and $I_{22}$, simultaneously. We will end up with the block matrix $\left( {\begin{array}{*{20}c}
S(z) & D(z)\\
0 & T(z)\\
\end{array} } \right)$ such that matrices $S(z)$, $D(z)$, and $T(z)$ satisfy (\ref{eq22}) and (\ref{eq24}).

Note that the proof of Lemma \ref{lemma36} below is an example of the diagonalization of a matrix polynomial by means of elementary transformations.

The following definition is in sync with definitions of root functions in \cite{GS, GLR}:

\begin{definition}\label{defition22} Let $L(z)$ be an $n\times n$ matrix polynomial. If an n-dimensional vector function $\varphi (z)$ satisfies 

\begin{enumerate}[(i)]
\item $\varphi (\alpha )\ne 0,$
\item $\left( L\left( z \right)\varphi \left( z \right) \right)^{\left( l \right)}\to 0$, when $z\to \alpha ,\, 0\le l\leq m-1$, 
\end{enumerate}
then we say that the function $\varphi (z)$ is a root function of $L(z)$ of order at least $m$ at the critical point (or zero) $\alpha $. If $m$ is the maximal number that satisfies (ii), then we say that $\varphi (z)$ is a root function of order $m$ at $\alpha$.
\end{definition}

The following characterization of root polynomials is a straightforward consequence of (\ref{eq120}); see the reasoning in \cite[p. 29]{GLR}.

\begin{proposition}\label{proposition23} Let $\alpha \in \sigma (L)$ be a zero of order $k$ of the polynomial $d_{i}(z)$, and $m \in \mathbb{N},\,1\leq m\leq k$. A function $\, \varphi (z)$ is a root function of 
$L(z)$ corresponding to $\alpha $ of order $m$ if and only if it is of the form 
\begin{equation}
\label{eq27}
\varphi \left( z \right)=\sum\limits_{j=0}^{m-1} {(z-\alpha )}^{j} \varphi
_{j}\left( \alpha \right)+\left( z-\alpha \right)^{m}\tilde{\varphi }\left( z \right),
\end{equation}
where $\varphi_{0}\left( \alpha \right),\, ... ,\, \varphi_{m-1}\left( \alpha \right)$ is a Jordan chain of $L$ at $\alpha$, not necessarily maximal, and 
$\tilde{\varphi }\left( z \right)$ is a function such that $\tilde{\varphi }\left( \alpha \right)$ is 
\textbf{not} the Jordan vector $\varphi_{m}$.
\end{proposition}

If $\varphi (z)$ is the root function of order $k\in \mathbb{N}$, i.e., $m=k$ in Definition \ref{defition22}, then $\varphi (z)$ is called a \textit{canonical} root function of $L\left( z \right)$ at $\alpha $. Obviously, the canonical root function is not uniquely determined, but according to Proposition \ref{proposition23}, the order $k$ is uniquely determined by the order of the corresponding canonical Jordan chain $\varphi_{0}\left( \alpha \right),\, ... ,\, \varphi_{k-1}\left( \alpha \right)$. The eigenvector $\varphi \left( \alpha \right)$ of $L\left( z \right)$ is also called a \textit{root vector} of $L(z)$ at $\alpha$ (see \cite{Lu2}).

The following version of Definition \ref{defition22} in the more general environment of 
meromorphic matrix functions introduces a couple of new terms. 

\begin{definition}\label{defition24} Given an $n\times n$ meromorphic matrix function 
$Q\left( z \right)$ and a holomorphic vector function $\psi \left( z 
\right)$ that satisfy

\begin{enumerate}[(i)]
\item $\varphi \left( z \right):=Q\left( z \right)\psi \left( z \right)\to c, c \ne 0, \, c \ne \infty ,\, if\, z\to \alpha ,$
\item ${\psi^{\left( l \right)} \left( z \right)}\to 0,\, if\, z\to \alpha ,\, 0\le l\leq m-1$, 
\end{enumerate}
then $\alpha \in \mathbb{C}$ is called a pole of $Q\left( z \right)$. The function 
$\varphi (z)$ is called the pole function of $Q(z)$ corresponding to 
$\alpha $, and $\psi \left( z \right)$ is called a pole cancellation 
function of $Q(z)$ of order at least $m$. If $m\in \mathbb{N}$ is the maximal number with property (ii) for the function $\psi$, then we say that functions $\varphi (z)$ and $\psi (z)$ are pole and pole cancellation functions, respectively,  of $Q\left( z \right)$ at $\alpha $ of order $m$. 
\end{definition}

For a definition of the pole cancellation function in the general case, when a generalized pole may be embedded in the singularities of an operator valued generalized Nevanlinna function $Q(z)$, refer to  \cite[Definitions 3.1 and 3.2]{BLu}. 

It is evident that $\varphi (z)$ is a root function of the invertible matrix polynomial 
$L(z)$ if and only if it is a pole function of $\hat{L}\left( z 
\right)=-{L(z)}^{-1}$ of the same order.
\\

\textbf{2.2.} According to (\ref{eq26}), for each $i \in \left\{ \mathrm{1,\, 
\mathellipsis ,\, }n \right\}$, the polynomial $d_{i}\left( z \right)$ might 
have multiple zeros. Also, a selected $\alpha \in \sigma (L)$ may be a zero of multiple polynomials $d_{i}\left(z \right) $ of the matrix $D(z)$ given by (\ref{eq24}). Therefore, for the diagonal matrix $D(z)$ we can introduce an index set as follows:
\[
\Omega \left( \alpha \right):=\left\{ i\in \mathbb{N}:\, d_{i}\left( \alpha \right)=0 \right\}.
\]
It is known that for any zero $\alpha $ of $d_{i}\left(z \right) $ of order $k$ there exists an eigenvector and a canonical Jordan chain of length $k$. In the following proposition, we will see how to find that eigenvector $\varphi_{0}\left( \alpha \right)$ and the corresponding canonical root polynomial $\varphi \left( z \right)$ by means of elementary transformations of columns of the matrix $L(z)$. Then we can easily obtain the entire canonical Jordan chain $\varphi_{0}\left( \alpha \right),\, ... ,\, \varphi_{k-1}\left( \alpha \right)$ by means of Proposition \ref{proposition23}. 

Let us denote the $i^{th}$ vector column of the matrix ${T\left( z 
\right)}$ by $T_{i}\left( z \right)$ and the $i^{th}$ vector 
column of the matrix $S^{-1}\left( z \right)$ by $\hat{S}_{i}\left( z \right)$. Then 
those two matrices can be represented as
\[
T\left( z \right)=\left[ T_{1}\left( z \right),\mathellipsis 
,\, T_{n}\left( z \right) \right]\, \wedge S^{-1}\left( z \right)=\left[ 
\hat{S}_{1}\left( z \right),\mathellipsis ,\, \hat{S}_{n}\left( z \right) \right].
\]

Now we have all necessary concepts to prove the following proposition. 

\begin{proposition}\label{proposition26} Let $\alpha \in \sigma \left( L \right)$, and let a diagonal form $D(z)$ of $L(z)$ be given by (\ref{eq24}). 

Let a scalar function $d_{i}\left( z \right)$ in representation (\ref{eq24}) be such that $ i \in \Omega \left( \alpha \right)$. Then the vector function 
\begin{equation}
\label{eq28}
\varphi_{i}\left( z \right):=T_{i}\left( z \right), \,  i \in \Omega \left( \alpha \right),
\end{equation}
is a canonical root polynomial of $L\left( z \right)$ corresponding to the 
eigenvalue $\alpha $ and the eigenvector $\varphi_{i,0}\left( 
\alpha \right)=T_{i}\left( 
\alpha \right)$. 

The vectors $\varphi_{i,0}\left( \alpha \right):=\varphi_{i}\left( \alpha \right)=T_{i}\left( \alpha \right)\ne 0$ are eigenvectors of $L(z)$ at $\alpha $ for all $i \in \Omega \left( \alpha \right)$. They are linearly independent, and 
\[
\ker {L\left( \alpha \right)}=span \left\{ T_{i}\left( \alpha 
\right):\, \, i \in \Omega \left( \alpha \right) \right\}.
\]
The canonical polynomials $\varphi_{i}\left( z \right),\, i \in \Omega \left( \alpha \right)$, are not necessarily uniquely determined; they depend on the selection of elementary transformations. The orders $k_{i}$ of $\varphi_{i}\left( z \right)$ are uniquely determined by the matrix $D(z)$. 
\end{proposition}

\textbf{Proof.} From (\ref{eq22}) it follows
\begin{equation}
\label{eq210}
L\left( z \right)T\left( z \right)={S\left( z \right)}^{-1}D\left( z 
\right),
\end{equation}
where ${S\left( z \right)}^{-1}$ is also a matrix polynomial. Indeed, because $S\left( z \right)$ is a matrix polynomial and the elements of ${S\left( z \right)}^{-1}$ are minors of ${S\left( z \right)}$ 
divided by the constant $\det\, S\left( z \right)$, ${S\left( z \right)}^{-1}$ is a matrix polynomial as well.

Now, let us select $\alpha \in \sigma (L)$ and introduce the n-dimensional vector $P_{i}\left( \alpha \right)$ with the $i^{th}$ component equal to 1, if $d_{i}\left( \alpha \right)=0$, and all other components equal to zero. Symbolically
\begin{equation}
\label{eq212}
P_{i}\left( \alpha \right)=\left[ p_{1}\left( \alpha 
\right),\mathellipsis ,\, \, p_{n}\left( \alpha \right) \right]^{T};\, 
p_{i}\left( \alpha \right)=\left\{ {\begin{array}{*{20}c}
1,\, if\, d_{i}\left( \alpha \right)=0. \\
0,\, othervise.\\
\end{array} } \right.
\end{equation}
Because matrices on the right and left sides of (\ref{eq210}) have identical 
elements, the corresponding columns will be 
identical too. Therefore, for the fixed $\alpha \in \sigma (L)$, and $i\in 
\left\{ 1,\, \mathellipsis ,\, n \right\}$ such that ${\, d}_{i}\left( 
\alpha \right)=0$,
\begin{equation}
\label{eq214}
L\left( z \right)T\left( z \right)P_{i}(\alpha )={S\left( z \right)}^{-1}D\left( z \right)P_{i}\left( \alpha \right), \forall i\in \Omega \left( \alpha \right). 
\end{equation}
By substituting (\ref{eq210}) and (\ref{eq212}) into the above equation, we get
\begin{equation}
\label{eq216}
L\left( z \right)T_{i}\left( z \right)=\hat{S}_{i}\left( z \right)\, 
d_{i}\left( z \right), \, \forall i \in \Omega \left( \alpha 
\right).
\end{equation}

We can now prove that the function $\varphi_{i}$ given by (\ref{eq28}) satisfies conditions (i) and (ii) of Definition \ref{defition22}.

(i) Because $\det\, T(z) \ne 0$ is a constant number, the vectors $T_{i}\left( z \right),\, i=1,\, \mathellipsis ,\, n$, are linearly independent vectors for every $z\in \mathbb{C}$. Therefore, $T_{i}\left( z \right) \ne 0,\, \forall z \in \mathbb{C}, i=1,\, \mathellipsis ,\, n$, and for the selected eigenvalue $\alpha$, the vectors 
\[
\varphi_{i}\left( \alpha \right)=T_{i}\left( \alpha \right)\ne 0,\, \, \forall i\in \Omega \left( \alpha \right).
\]
This proves that Definition \ref{defition22} (i) is satisfied.

(ii) Let us now observe representations (\ref{eq26}). Let us arbitrarily select an index $i\in \lbrace 1, 2, ..., n \rbrace$ and a zero $z=\alpha $ of the polynomial $d_{i}\left( z\right)$. Let us denote the multiplicity of the zero $z=\alpha $ by $k_{i}$. Then $d_{i}\left( z \right)=\left( z-\alpha \right)^{k_{i}}\tilde{d}_{i}\left( z \right)$, where $\tilde{d}_{i}\left( \alpha \right)\ne 0$. This, (\ref{eq28}), and (\ref{eq216}) imply
\begin{equation}
\label{eq218}
L\left( z \right)\varphi_{i}\left( z \right)=\left( z-\alpha 
\right)^{k_{i}}\hat{S}_{i}\left( z \right)\tilde{d}_{i}\left( z \right).
\end{equation}
Since $\hat{S}_{i}\left( z \right)$ and $\tilde{d}_{i}\left( z \right)$ are holomorphic at $z=\alpha$ we have
\[
\left( L\left( z \right)\varphi_{i}\left( z \right) \right)^{\left( l 
\right)}\to 0,\, if\, z\to \alpha ,\, 0\le l\le k_{i}-1,
\]
Therefore, $\varphi_{i}\left( z \right)=T_{i}\left( z \right)$ satisfies Definition \ref{defition22} (ii), i.e, it is a root function of order at least $k_{i}$. 

It remains to see that $k_{i}$ is the maximal number that satisfies Definition \ref{defition22} (ii). 
Indeed, let us assume that it holds
\[
\left( L\left( z \right)\varphi_{i}\left( z \right) \right)^{\left( k_{i} 
\right)} \to 0,\, \left( z\to \alpha \right).
\]
If we differentiate the right-hand side of (\ref{eq218}) $k_{i}$ times, we will get:
\[
\left( \left( z-\alpha \right)^{k_{i}}\hat{S}_{i}\left( z \right)\tilde{d}_{i}\left( z \right)\right)^{( k_{i})}=\left( z-\alpha \right)O(z)+ \hat{S}_{i}\left( z \right)\tilde{d}_{i}\left( z \right),
\]
where $O(z)$ is a holomorphic function. This means
\[
\hat{S}_{i}\left(z\right)\tilde{d}_{i}\left( z\right) \to \hat{S}_{i}\left( \alpha \right)\tilde{d}_{i}\left( \alpha \right)=0,\, \left( z\to \alpha \right).
\]
Because $\tilde{d}_{i}\left( \alpha \right) \neq 0 $, it follows $\hat{S}_{i}\left( \alpha \right)=0$. This 
is in contradiction with the assumption $det\, \tilde{S}\left( z \right)=\left( det\, S \left( z \right)\right)^{-1}=c\ne 0$. This proves that Definition \ref{defition22} is satisfied with $k_{i}$ as the maximal number.

The vectors $\varphi_{i}\left( \alpha \right), i \in \Omega(\alpha)$, are linearly independent eigenvectors of $L(z)$ at $\alpha$ as columns of the regular matrix $T(\alpha)$, and by definition of $\Omega (\alpha)$ they span entire $\ker L(\alpha)$. 

The final statement of the proposition follows from the fact that the matrix $T(z)$ generally depends on the selection of the elementary transformations that lead to the particular diagonalization $D(z)$, while orders $k_{i}$ are uniquely determined for the matrix $D(z)$ because those are the orders of the solutions of the equations $d_{i}(z)=0$.   \hfill $\square$

\begin{corollary}\label{corollary210} Let $D\left( z \right)$ be a diagonal form (\ref{eq24}) of the matrix polynomial $L\left( z \right)$ and let $S \left( z \right)$ and $T \left( z \right)$ be matrices of elementary transformations that satisfy (\ref{eq22}). Then for every $i=1, ...,n $ the function $\varphi_{i} \left( z \right):=T_{i}\left( z \right) $ is the root function of $L\left( z \right) $ at all eigenvalues that satisfy $d_{i}\left( z \right)=0$. 
\end{corollary}

The advantage of the method of finding root functions for $L(z)$ presented in the above proposition and the corollary is that it is based on a sequence of elementary transformations of matrices, a topic typically covered at the college level of Linear Algebra. Additionally, it has several other benefits. For instance, we do not need to solve the algebraic equation $\det {L(z)}=0$, which is of degree $nl$, to find eigenvalues of $L$. We do not even need to find $\det {L(z)}$. Instead, once we obtain the diagonal matrix $D(z)$ we only have to solve algebraic equations $d_i(z)=0$, $i=1,\ldots,n$, which are of smaller orders. Furthermore, although the proof uses properties of $\det {S(z)}$, $\det {T(z)}$, and $S(z)^{-1}$, we do not need to calculate any of these objects.

In the following corollary, we determine the set of canonical root functions without solving any of the equations $d_i(z)=0$.

\textbf{2.3.} Let us introduce the function 
\begin{equation}
\label{eq220}
\psi_{i}\left( z \right):=L\left( z \right)\varphi_{i}\left( z \right)=\left( z-\alpha 
\right)^{k_{i}}\hat{S}_{i}\left( z \right)\tilde{d}_{i}\left( z \right),
\end{equation}
where, as we know, $\hat{S}_{i}\left( \alpha \right)\tilde{d}_{i}\left( \alpha 
\right)\ne 0.$ This means that the eigenvalue $\alpha $ of $L\left( z 
\right)$ is a zero of $\psi_{i}\left( z \right)$ of exact order $k_{i}$. Then it holds
\[
\hat{L}\left( z \right){\psi_{i}\left( z \right)} \rightarrow -\varphi_{i}\left( \alpha \right)\ne 0 \ as \ z\to \alpha.
\]
Hence, the function $\psi_{i}\left( z \right)$ defined by (\ref{eq220}) is a \textit{pole cancellation function} of $\hat{L}\left( z \right)$ of order $k_{i}$ that corresponds to the canonical root function $\varphi_{i}\left( z \right)$. 

In Example \ref{example214}, we will show how to solve the homogeneous system of linear differential equations (\ref{eq18}) by means of Proposition \ref{proposition26}. 

\begin{example}\label{example214} Solve  the homogeneous system  of differential equations 
\begin{equation}
\label{eq222}
\left( {\begin{array}{*{20}c}
0 & \mathrm{\, }\frac{d}{dt} & 0\\
1 & 0 & \mathrm{\, }\frac{d}{dt}\\
0 & 0 & \mathrm{\, }\frac{d^{3}}{{dt}^{3}}-\frac{d^{2}}{{dt}^{2}}\\
\end{array} } \right)\left( {\begin{array}{*{20}c}
u_{1}\left( t \right)\\
u_{2}\left( t \right)\\
u_{3}\left( t \right)\\
\end{array} } \right)=\left( {\begin{array}{*{20}c}
0\\
0\\
0\\
\end{array} } \right).
\end{equation}
\begin{enumerate}[(i)]
\item Find the canonical systems of root function and Jordan chains of the matrix polynomial $L\left( z \right)$ corresponding to (\ref{eq222}) by means of Proposition \ref{proposition26}. 
\item Find the general solution of the system (\ref{eq222}) using the results of (i).
\end{enumerate}
\end{example} 

(i) Clearly
\[
L\left( z \right)=\left( {\begin{array}{*{20}c}
0 & z & 0\\
1 & 0 & z\\
0 & 0 & z^{3}-z^{2}\\
\end{array} } \right).
\]
In this example, the identity (\ref{eq22}) is
\[
\left( {\begin{array}{*{20}c}
0 & 1 & 0\\
1 & 0 & 0\\
0 & 0 & 1\\
\end{array} } \right)\left( {\begin{array}{*{20}c}
0 & z & 0\\
1 & 0 & z\\
0 & 0 & z^{3}-z^{2}\\
\end{array} } \right)\left( {\begin{array}{*{20}c}
1 & 0 & -z\\
0 & 1 & 0\\
0 & 0 & 1\\
\end{array} } \right)=\left( {\begin{array}{*{20}c}
1 & 0 & 0\\
0 & z & 0\\
0 & 0 & z^{3}-z^{2}\\
\end{array} } \right)=D\left( z \right).
\]
According to (\ref{eq210}) we get
\[
\left( {\begin{array}{*{20}c}
0 & z & 0\\
1 & 0 & z\\
0 & 0 & z^{3}-z^{2}\\
\end{array} } \right)\left( {\begin{array}{*{20}c}
1 & 0 & -z\\
0 & 1 & 0\\
0 & 0 & 1\\
\end{array} } \right)=\left( {\begin{array}{*{20}c}
0 & 1 & 0\\
1 & 0 & 0\\
0 & 0 & 1\\
\end{array} } \right)\left( {\begin{array}{*{20}c}
1 & 0 & 0\\
0 & z & 0\\
0 & 0 & z^{3}-z^{2}\\
\end{array} } \right).
\]
First, let us find the root polynomials at $\alpha 
=0$. In that case $\Omega \left( 0 \right):=\left\{ 2,\, 3 \right\}$. 

For $i=2$, according to formula (\ref{eq28}), we have
\[
P_{2}\left( 0 \right)=\left( {\begin{array}{*{20}c}
0\\
1\\
0\\
\end{array} } \right).
\]
According to (\ref{eq214}) it follows
\begin{equation}
\label{eq224}
\left( {\begin{array}{*{20}c}
0 & z & 0\\
1 & 0 & z\\
0 & 0 & z^{3}-z^{2}\\
\end{array} } \right)\left( {\begin{array}{*{20}c}
0\\
1\\
0\\
\end{array} } \right)=\left( {\begin{array}{*{20}c}
z\\
0\\
0\\
\end{array} } \right).
\end{equation}
According to definition (\ref{eq220}) of $\psi_{i}(z)$, $\psi_{2}\left( z \right)=\left( 
{\begin{array}{*{20}c}
z\\
0\\
0\\
\end{array} } \right)$ and $\varphi_{2}\left( z \right)=\left( {\begin{array}{*{20}c}
0\\
1\\
0\\
\end{array} } \right)$. Because $z=0$ is a zero of the first order of $\psi_{2}\left( z \right)$, we know that the Jordan chain of the corresponding eigenvector will consist only of that eigenvector $\varphi_{2,0}$. Therefore $\varphi_{2}\left( z \right)=\varphi_{2,0}=\left( 
{\begin{array}{*{20}c}
0\\
1\\
0\\
\end{array} } \right)$ and $\tilde{\varphi }_{2}\left( z \right)=\left( 
{\begin{array}{*{20}c}
0\\
0\\
0\\
\end{array} } \right)$. This is consistent with the Proposition \ref{proposition23} because 
$\tilde{\varphi }_{2}\left( 0 \right)$ is not a Jordan vector.

Next, let us find the root polynomial of $L(z)$ at the eigenvalue $\alpha =0$, and 
$i=3$. Then $P_{3}\left( 0 \right)=\left( {\begin{array}{*{20}c}
0\\
0\\
1\\
\end{array} } \right)$.

According to (\ref{eq214}) it follows
\begin{equation}
\label{eq226}
\left( {\begin{array}{*{20}c}
0 & z & 0\\
1 & 0 & z\\
0 & 0 & z^{3}-z^{2}\\
\end{array} } \right)\left( {\begin{array}{*{20}c}
-z\\
0\\
1\\
\end{array} } \right)=\left( {\begin{array}{*{20}c}
0\\
0\\
z^{3}-z^{2}\\
\end{array} } \right).
\end{equation}
Because $\alpha =0$ is the zero of the second order of $\psi_{3}\left( z 
\right)=\left( {\begin{array}{*{20}c}
0\\
0\\
z^{3}-z^{2}\\
\end{array} } \right)$, we know that there exists a Jordan chain of 
length $k=2$ at $\alpha =0$. When we express $\varphi_{3}\left( z \right)$ in 
exponents of $z$, we get
\[
\varphi_{3}\left( z \right)=\left( {\begin{array}{*{20}c}
-z\\
0\\
1\\
\end{array} } \right)=\left( {\begin{array}{*{20}c}
0\\
0\\
1\\
\end{array} } \right)+\left( {\begin{array}{*{20}c}
-1\\
0\\
0\\
\end{array} } \right)z.
\]
Therefore, the Jordan chain is: $\varphi_{3,0}=\left( 
{\begin{array}{*{20}c}
0\\
0\\
1\\
\end{array} } \right), \, \varphi_{3,1}=\left( {\begin{array}{*{20}c}
-1\\
0\\
0\\
\end{array} } \right)$. Then $\tilde{\varphi }_{3}\left( z \right)=\left( 
{\begin{array}{*{20}c}
0\\
0\\
0\\
\end{array} } \right)$.

Finally, let us find the root polynomial $\varphi \left( z \right)$ at $\alpha 
=1.$ Again, we use $P_{3}\left( 0 \right)=\left( {\begin{array}{*{20}c}
0\\
0\\
1\\
\end{array} } \right)$. Therefore we obtain again (\ref{eq226}) and 
$\varphi_{3}\left( z \right)=\left( {\begin{array}{*{20}c}
-z\\
0\\
1\\
\end{array} } \right)$. When we express $\varphi_{3}\left( z \right)$ at $\alpha=1$ in 
exponents of $z-1$, we get
\[
\varphi_{3}\left( z \right)=\left( {\begin{array}{*{20}c}
-z\\
0\\
1\\
\end{array} } \right)=\left( {\begin{array}{*{20}c}
-1\\
0\\
1\\
\end{array} } \right)+\left( {\begin{array}{*{20}c}
-1\\
0\\
0\\
\end{array} } \right)\left( z-1 \right).
\] 
Because $\alpha =1$, is a zero of the first order of $\psi_{3}\left( z 
\right)$, we conclude that the Jordan chain consists only of $\varphi 
_{3,0}=\left( {\begin{array}{*{20}c}
-1\\
0\\
1\\
\end{array} } \right)$. Then $\tilde{\varphi }_{3}\left( z \right)=\left( 
{\begin{array}{*{20}c}
-1\\
0\\
0\\
\end{array} } \right)$ and it is not a Jordan vector. 

Note, the function $\varphi_{3}\left( z \right)=T_{3}\left( z \right)$ is the root polynomial at two eigenvalues of $L(z)$, i.e., at both zeros of the polynomial $d_{3}\left( z \right)=z^{3}-z^{2}$;  $\varphi_{3}\left( z \right)$ is a canonical root function of $L(z)$ at $z=0$ of order $k=2$, and a canonical root function of $L(z)$ at $z=1$ of order $k=1$. This is consistent with Corollary \ref{corollary210}.

(ii) We have three linearly independent solutions. They are given by (\ref{eq118}). For $\alpha=0$ we have the solution $u_{0,1}\left( t 
\right)=\varphi_{0}=\left( {\begin{array}{*{20}c}
0\\
1\\
0\\
\end{array} } \right)$. The second solution we have again for $\alpha =0$:
\[
u_{0,2}\left( t \right)=t\varphi 
_{0}+\varphi_{1}=\left( {\begin{array}{*{20}c}
0\\
0\\
1\\
\end{array} } \right)t+\left( {\begin{array}{*{20}c}
-1\\
0\\
0\\
\end{array} } \right)=\left( {\begin{array}{*{20}c}
-1\\
0\\
t\\
\end{array} } \right).
\]
For $\alpha =1$, we have the solution: $u_{1,3}\left( t \right)=\left( 
{\begin{array}{*{20}c}
-1\\
0\\
1\\
\end{array} } \right)e^{t}$. Then the general solution is $u(t)=\left( {\begin{array}{*{20}c}
u_{1}\left( t \right)\\
u_{2}\left( t \right)\\
u_{3}\left( t \right)\\
\end{array} } \right)=C_{1}u_{0,1}\left( t \right)+ C_{2}u_{0,2}\left( t \right)+C_{3}u_{1,3}\left( t \right)$, where $ C_{i}, \, 1=1,2,3$, are complex constants.
\hfill $\square$

\section{Hermitian polynomials $L\left( z \right)$}\label{s6}
\textbf{3.1.} Consider the \textit{Hermitian matrix polynomial} defined by
\begin{equation}
\label{eq32}
L\left( z \right):=A_{l}z^{l}+\mathellipsis +A_{1}z+A_{0},
\end{equation}
where $A_{i}, i=0.1,..., l$ are $n\times n$ Hermitian matrices, and $l=\deg\, L(z)\, \, \ge 1$. Recall that an 
$n \times n$-complex matrix $A=\left( a_{ij} \right)$ is \textit{self-adjoint (Hermitian}) if it satisfies 
$a_{ij}=\overline{a_{ji}}$.  If $\det {L\left( z \right)}\not\equiv 0$, then the polynomial is invertible and we use the notation $\hat{L}\left( z \right)=-{L\left( z \right)}^{-1}$.

It is well-known and easy to see that the functions $L\left( z \right)$ and $\hat{L}\left( z \right)$ are \textit{matrix generalized Nevanlinna functions} with the same negative index $\kappa \in \mathbb{N}\cup \left\{ 0 \right\}$, symbolically denoted as $L,\, \hat{L}\in N_{\kappa }^{n\times n}$. 

Since the matrix polynomials are crucial for solving systems of linear ODEs, our focus in this paper is on the investigating of the properties of $L(z)$ and $\hat{L}(z)$. However, some of the following statements will be proven in the more general context of generalized Nevanlinna matrix functions $N_{\kappa }^{n\times n}$, or operator-valued functions $ N_{\kappa }(\mathcal{H})$.

Recall that the function $Q\in N_{\kappa }(\mathcal{H}) $ is said to be convergent at $\infty$ if the limit 
\begin{equation}
\label{eq33}
\lim\limits_{z\to \infty}{{Q\left( z \right)}}=S \in \mathcal{L}(\mathcal{H}).
\end{equation}
exists, and it is called holomorphic at $\infty$ if
\begin{equation}
\label{eq36}
Q'(\infty):=\lim\limits_{z\to \infty}{zQ\left( z \right)}
\end{equation}
exists in the Banach space of bounded operators $ N_{\kappa }(\mathcal{H}) $. We will need the following representation from \cite{B2}.
\begin{lemma}\label{lemma32} \cite[Lemma 3]{B2}.
A function $Q\in N_{\kappa }(\mathcal{H}) $ is holomorphic at $\infty$ if and only if $Q(z)$ has minimal 
representation 
\begin{equation}
\label{eq34}
Q\left( z \right)=\Gamma^{+}\left( A-z \right)^{-1}\Gamma , \,z \in \dom\, Q ,
\end{equation}
with a bounded self-adjoint operator $A$ in a Ponteyagin space $\mathcal{K}$. In this case $Q'(\infty)=-\Gamma^{+}\Gamma$.
\end{lemma}
\begin{proposition}\label{proposition32.1} A function $Q\in N_{\kappa }(\mathcal{H}) $ is convergent at $\infty$ and the function
\[
\tilde{Q}(z):= Q(z)-S
\]
is holomorphic at $\infty$ if and only if $Q(z)$ has the representation
\begin{equation}
\label{eq37}
Q\left( z \right)=S+\Gamma^{+}\left( A-z \right)^{-1}\Gamma ,\,  {S=S}^{\ast }\in \mathcal{L} \left( \mathcal{H} \right) 
\end{equation}
where $A$ is a self-adjoint bounded operator in the minimal Pontryagin space $\left( \mathcal{K},\left[ .,.\right] \right) $, $\Gamma :\mathcal{H}\to \mathcal{K}$ and $\Gamma^{+}:\mathcal{K}\to \mathcal{H}$ are linear 
operators. 
\end{proposition} 
\textbf{Proof.} If the function $Q$ satisfies (\ref{eq33}) with $S=0$, then $Q(z)= \tilde{Q}(z)$ is holomorphic at $\infty$ and (\ref{eq37}) with $S=0$ follows from Lemma \ref{lemma32}. If $S \neq 0$, then according to assumption $\tilde{Q}$ is holomorphic at $\infty$. Then, Lemma \ref{lemma32} implies $Q(z)-S=\Gamma^{+}\left( A-z \right)^{-1}\Gamma $, which proves  \eqref{eq37}. 

Conversely, assume that representation \eqref{eq37} holds with bounded self-adjoint operator $A$ in the Pontryagin space $\mathcal{K}$. Then, according to Lemma \ref{lemma32}, $\tilde{Q}$ is holomorphic at $\infty$, with
\[  
\tilde{Q}^{'}\left( \infty \right)=-\Gamma^{+}\Gamma.
\] 
Then obviously the function $Q$ is convergent at $\infty$ with limit $S$. \hfill $\square$ 

\begin{corollary}\label{corollary32.2} Let $L\left( z \right)$, $l\geq 1$, be a scalar polynomial with real coefficients. Then $\hat{L}(z)$ has a minimal representation (\ref{eq34}), i.e.,
\[
\hat{L}(z)=\left[ v, \left( A-z \right)^{-1}v\right], 
\]
where $v \in \dom\, A=\mathcal{K}$ is a fixed element that satisfies
\[
\mathcal{K}=c.l.s.\left\{( A-z)^{-1}v :z\in \dom\,Q \right\}.
\]
\end{corollary}
This corollary directly follows from Lemma \ref{lemma32} for $\mathcal{H}=\mathbb{C}$ and $Q=\hat{L}$. For the corresponding statement in the general class of scalar generalized Nevanlinna functions see \cite[Satz 1.5]{KL1}.
\\

\textbf{3.2.} Corollary \ref{corollary32.2} cannot be generalized to Hermitian matrix polynomials $L(z)$ without an additional condition. In this section, in Proposition \ref{proposition34}, we will identify the condition under which we can prove the corresponding statement for Hermitian matrix polynomial $L(z)$. 

\begin{lemma}\label{lemma33} Let $L(z)$ be an $n \times n$ matrix polynomial given by (\ref{eq32}), not necessarily Hermitian, with characteristic polynomial $\chi_{n}(z):=\det L(z) \not\equiv 0$. Then:
\begin{enumerate}[(i)]
\item $\deg \chi_{n}(z) \leq nl$.
\item If the polynomial $L$ is monic, then $L(z)$ is invertible and $\deg \chi_{n}(z)=nl$.
\end{enumerate} 
\end{lemma}

\textbf{Proof.} By definition of the inverse matrix ${L\left( z \right)}^{-1}$, every element of $\hat{L}(z)= -{L\left( z \right)}^{-1}$ is the ratio of two polynomials, i.e.,
\begin{equation}
\label{eq37.2}
\hat{L}_{ij}(z)=-\frac{(-1)^{i+j}m_{ji}(z)}{\chi_{n} (z)},\, i.j=1,...,n,
\end{equation}
where $m_{ji}(z), \, j,i=1, ..., n$, denotes the determinant of the corresponding minor of $L(z)$ of order $n-1$.

(i) We will use mathematical induction on $n\in \mathbb{N}$, for an arbitrarily selected $l\in \mathbb{N}$. The statement $\deg \chi_{n}(z) \leq nl$ obviously holds for $n=1, 2$. Let us assume that for matrix polynomials of order $n-1$ and degree $l\geq 1$ the characteristic polynomials $\chi_{n-1}(z)$ satisfy
\[
\deg \chi_{n-1}(z) \leq (n-1)l.
\]
The polynomial $\chi_{n}(z)=\det L(z)$ is by definition a sum of $n$ products of the form $(-1)^{i+j}L_{ij}(z)m_{ij}(z)$, where $L_{ij}(z)$ are elements of $L(z)$ and $m_{ij}(z)$ are determinants of the corresponding minors of order $n-1$. The elements $L_{ij}(z)$  are polynomials of degrees $ k =0, 1, ..., l$, while $\deg m_{ij}(z) \leq (n-1)l$, according to the assumption of induction. Therefore, for the degree of each summand of the determinant $\chi_{n}(z)$ we have $k +(n-1)l \leq nl$. This proves that $\deg \chi_{n}(z)\leq nl$ for all $n=1,2,...$. 

(ii) By definition of the monic polynomial $L(z)$, the exponents $z^{l}$ appear only on the main diagonal, while all other elements of $L(z)$ are polynomials with degrees less than $l$. Let $\left(  \sigma(1),\sigma(2), ... ,\sigma(n) \right) $ denote a permutation of the set $S_{n}:=\lbrace1,2, ... ,n \rbrace $. According to Leibniz formula, $\det L(z)$ is a sum of products of the form  $\pm L(z)_{1\sigma (1)}L(z)_{2\sigma (2)}...L(z)_{n\sigma (n)} $, where the indexes $\left(  \sigma(1),\sigma(2), ... ,\sigma(n) \right) $ run over all $n!$ permutations of $S_{n}$. Therefore, all $n!$ summands are products of exactly $n$ various polynomials $L_{ij}(z)$, each with degree less than $l$, except the elements $L(z)_{ii}=z^{l}, i=1, 2, ..., n$. Thus, $\det L(z)$ will contain only one summand $L(z)_{11}L(z)_{22}...L(z)_{nn}= z^{nl}$, the product of $n$ elements from the main diagonal. This means $\deg \chi_{n}(z)= nl$. This also means $\det L(z) \not\equiv 0 $, which proves the invertibility of $L(z)$.  \hfill $\square$

According to Lemma \ref{lemma33} the condition
\begin{equation}
\label{eq38}
\frac{\deg m_{ij}(z)}{\deg \chi_{n} (z)}<1,\, \forall i.j=1,...,n,
\end{equation}
holds whenever $\deg \chi_{n} (z)=nl$. The following example proves that the condition (\ref{eq38}) need not to hold when $\deg \chi_{n} (z)\neq nl$. Moreover, it is an example of the inverible matrix polynomial $L(z)$ for which the inverse $\hat{L}(z)$ is not convergent at $\infty$, and therefore not holomorphic at $\infty$.

\begin{example}\label{example35} Given the Hermitian polynomial $L(z)= \left( {\begin{array}{*{20}c}
z^{3} & z\\
z & 0\\
\end{array} } \right) $, we have $\chi_{2}(z)=-z^{2}.$ Therefore, $\deg \chi_{2}(z)=2 < nl=6$. Moreover, $\deg \chi_{2}(z)=2 < l=3$. A very important consequence of this property is that the inverse $\hat{L}(z):=-L^{-1}(z)= \left( {\begin{array}{*{20}c}
0 & z^{-1}\\
z^{-1} & z\\
\end{array} } \right) $ is not a function convergent at $\infty$, hence, it does not satisfy condition (\ref{eq38}). 

Additionally, this is an example where the representing relation $A$ of $\hat{L}$ has an eigenvalue at $\infty$, as the function $\hat{L}(-\varsigma^{-1})$ has an eigenvalue at $\varsigma=0$. Therefore, the function $\hat{L}$ can not have the simplified representation \eqref{eq37}, but only the most general representation \eqref{eq12}.
\end{example}
\begin{proposition}\label{proposition34} Let $L\left( z \right)$, $l\geq 1$, be an invertible Hermitian 
$n\times n. \, n\geq 2$ matrix polynomial. The condition
\begin{equation}
\label{eq310}
\frac{\deg m_{ij}(z)}{\deg \chi_{n} (z)}\leq1,\, \forall i.j=1,...,n,
\end{equation} 
is satisfied if and only if $\hat{L}\left( z \right)$ has the representation 
\begin{equation}
\label{eq312}
\hat{L}\left( z \right)=S+\Gamma^{+}\left( A-z \right)^{-1}\Gamma ,\,  {S=S}^{\ast }\in \mathcal{L} \left( \mathbb{C}^{n} \right) 
\end{equation}
where $A$ is a self-adjoint bounded operator in the minimal Pontryagin space $\left( \mathcal{K},\left[ .,.\right] \right) $, and $\Gamma :\mathbb{C}^{n}\to \mathcal{K}$ and $\Gamma^{+}:\mathcal{K}\to \mathbb{C}^{n}$ are linear 
operators. 
\end{proposition} 

\textbf{Proof.} Assume that the inverse matrix function $\hat{L}$ satisfies condition (\ref{eq310}). Since the elements of the matrix $\hat{L}(z)$ are of the form \eqref{eq37.2}, with the degrees of the polynomials in the numerators being less then or equal to the degree of the polynomial in the denominator, the function $\hat{L}(z)$ converges at infinity to some matrix $S$. Consequently, the function $\tilde{L}(z):=\hat{L}(z)-S$ is holomorphic at $\infty$. According to Proposition \ref{proposition32.1}, the representation \eqref{eq312} holds.

Conversely, if the representation \eqref{eq312} holds, then the function $\tilde{L}(z):=\hat{L}(z)-S$ is holomorphic at $\infty$. Therefore it satisfies condition \eqref{eq38}, which implies that $\hat{L}$ satisfies condition \eqref{eq310}. \hfill $\square$


\textbf{3.3.} If $\alpha$ is a generalized pole of $Q\in N_{\kappa }(\mathcal{H}) $, then there exist eigenvectors and Jordan chains of rhe representing relation $A$ at $\alpha $, (see e.g. \cite[Definition 4.2]{Lu3}). The linear span of all Jordan vectors of $A$ at $\alpha $ is called \textit{algebraic eigenspace of} $A$ at $\alpha $ and it is denoted by $S_{\alpha }\left( A \right)$ (see e.g., \cite{IKL}). Obviously, $A\left( S_{\alpha }\left( A \right) \right)\subseteq S_{\alpha }\left( A \right)$. 

We already mentioned that eigenvalues of $L(z)$ and eigenvalues of the representing relation $A$ of $\hat{L}(z)$ in the  Krein-Langer type of representation coincide. Therefore, it is possible to investigate properties of $L(z)$ by investigating properties of the operator $A$, and vice versa. Unfortunately, it is very difficult to find an actual operator representation, i.e., the Pontraygin space $\mathcal{K}$, the self-adjoint relation $A$ in $\mathcal{K}$, and the operator $\Gamma:\mathcal{H} \rightarrow\mathcal{K}$ for a given generalized Nevanlinna function $Q\in N_{\kappa }(\mathcal{H}) $. The constructions in the cited papers are abstract and not easy to apply. The following Theorem \ref{theorem38} is one of the first practical results in that direction. In that theorem we will give a method how to obtain the representing operators $A$ and $\Gamma$ from the given meromorphic function $Q\in N_{\kappa }^{n \times n} $ that satisfies conditions of Proposition \ref{proposition32.1}. 

We will need the following lemma. 

\begin{lemma}\label{lemma36} Let $G=\left( {\begin{array}{*{20}c}
0 & \cdots & 1\\
\vdots & \ddots & \vdots \\
1 & \cdots & 0\\
\end{array} } \right)$ be a $k\times k$ matrix, i.e., $g_{i\left( k+1-i 
\right)}=1,\, i=1,\, \mathellipsis ,\, k$, all others $g_{ij}=0$.

\begin{enumerate}[(i)]
\item If $k=2l+1, l\in \lbrace0\rbrace\cup\mathbb{N}$, then the matrix $G$ has $l$ negative and $l+1$ positive eigenvalues, and the matrix $-G$ has $l$ positive and $l+1$ negative eigenvalues. 
\item If $k=2l$, then both matrices, $G$ and $-G$, have $l$ negative and $l$ positive eigenvalues.
\end{enumerate}
\end{lemma}

\textbf{Proof.} 

(i) $k=2l+1$. We will describe how to obtain a convenient diagonal matrix $B\left( z 
\right)$ equivalent with the matrix polynomial $G-zI$ by listing the first 
three steps, i.e., the first three elementary transformations. 

\begin{enumerate}[1)]
\item Reverse the order of rows or columns:

\[
G-zI=\left( {\begin{array}{*{20}c}
-z & \cdots & 1\\
\vdots & \ddots & \vdots \\
1 & \cdots & -z\\
\end{array} } \right)\to \left( {\begin{array}{*{20}c}
1 & \cdots & -z\\
\vdots & \ddots & \vdots \\
-z & \cdots & 1\\
\end{array} } \right).
\]
\item Multiply the first column by $z$ and add to the last:
\[
\left( {\begin{array}{*{20}c}
1 & \cdots & -z\\
\vdots & \ddots & \vdots \\
-z & \cdots & 1\\
\end{array} } \right)\to \left( {\begin{array}{*{20}c}
1 & \cdots & 0\\
\vdots & \ddots & \vdots \\
-z & \cdots & 1-z^{2}\\
\end{array} } \right).
\]
\item Multiply the first row by $z$ and add to the last:
\[
\left( {\begin{array}{*{20}c}
1 & \cdots & 0\\
\vdots & \ddots & \vdots \\
-z & \cdots & 1-z^{2}\\
\end{array} } \right)\to \left( {\begin{array}{*{20}c}
1 & \cdots & 0\\
\vdots & \ddots & \vdots \\
0 & \cdots & 1-z^{2}\\
\end{array} } \right).
\]
\item By repeating steps 2 and 3 with the second column and the second row, and so on, we obtain the diagonal matrix $B$ with
\[
b_{11}=\mathellipsis =b_{ll}=1,\, b_{\left( l+1 \right)\left( l+1 \right)}=1-z,
\]
\[
\, b_{\left( l+2 \right)\left( l+2 \right)}=1-z^{2},\, \mathellipsis ,\, b_{\left( 2l+1 \right)\left( 2l+1 \right)}=1-z^{2},
\]
i.e., 
\begin{equation}
\label{eq315}
B\left( z \right)=diag \, \left[ 1,\, \mathellipsis ,\, 1,\, \left( 1-z \right),\, \left( 1-z^{2} \right),\mathellipsis ,\left( 1-z^{2} \right) \right].
\end{equation}
Therefore, 
\[
B\left( z \right)=S\left( z \right)\left( G-zI \right)T\left( z \right),
\]
\end{enumerate}
where $S\left( z \right)$ is the product of elementary row transformations and $T \left( z \right)$ is the product of elementary column transformations. We know that the determinant of the product of 
matrices is equal to the product of the determinants. Therefore, $\det {S\left( 
z \right)}$ and $\det {T \left( z \right)}$ are constants. This further 
means that the matrices $B \left( z \right)$ and $G-zI$ have the same 
characteristic polynomials and eigenvalues. Now statement (i) follows 
from the above diagonal representation $B\left( z \right)$ of $G-zI$.

The statement (ii), when $k=2l$, follows by applying the above proof on each of the matrices $G$ and $-G$. \hfill $\square$

The following theorem, combined with Proposition \ref{proposition26}, will enable us to find a Pontryagin space and the representing operators $A$ and $\Gamma$ of $\hat{L}\left( z \right)$ for the representation (\ref{eq312}).

\begin{theorem}\label{theorem38} Let $Q\in N_{\kappa }^{n \times n} $ satisfies conditions of Proposition \ref{proposition32.1}, and let $A$ be the self-adjoint operator in the Pontryagin space $\left( \mathcal{K},\left[ .,.\right] \right) $ in the minimal representation (\ref{eq37}) of $Q$. 
\begin{enumerate}[(i)]
\item Let $\varphi \left( z \right)$ be a canonical root function of $\hat{Q}:=-Q^{-1}$ at 
$\alpha \in \sigma \left( Q \right)\cap \mathbb{R}$ of order $k$. Then there exists a Jordan chain of the operator $A$ at $\alpha $
\begin{equation}
\label{eq320}
f_{0},\, f_{1},\, \mathellipsis ,\, f_{k-1},
\end{equation}
such that it holds: 
\begin{equation}
\label{eq316}
\lim\limits_{z\to \alpha}{\frac{\hat{Q}\left( z \right)\varphi \left( z \right)}{{(z-\alpha )}^{j}}}=0 \,\left( j < k \right),
\end{equation}
\begin{equation}
\label{eq318}
\lim\limits_{z\to \alpha}{\frac{\left( \hat{Q}\left( z \right)\varphi \left( z \right),\varphi \left( z \right) \right)}{{(z-\alpha )}^{k}}=\left[ f_{0},f_{k-1} \right]}\ne 0.
\end{equation}
The space
\begin{equation}
\label{eq321}
S_{\varphi}:=l.s.\, \left\{ f_{0},\, f_{1},\, \mathellipsis ,\, f_{k-1} \right\}
\end{equation}
is a Pontryagin subspace of $K$, and the operator $A_{\varphi}:=A_{\vert S_{\varphi}}$, with Jordan block form 
\begin{equation}
\label{eq321.1}
A_{\varphi}=\left( {\begin{array}{*{20}c}
\alpha & 1 & ... & 0 & 0\\
0 & \alpha  & ... & 0 & 0\\
\vdots & \vdots & \vdots &\vdots & \vdots\\
0 & 0 & ... & \alpha & 1\\
0 & 0 & ... & 0 & \alpha\\
\end{array} } \right),
\end{equation}
is a symmetric operator with respect to the indefinite scalar product $\left[ .,. \right]_{S_{\varphi}}:=\left[ .,. \right]_{\vert S_{\varphi}}=\left( J_{\varphi}.,.\right)  $.
The fundamental symmetry of the scalar product $\left[ .,. \right]_{S_{\varphi}}$ is 
\begin{equation}
\label{eq321.15}
J_{\varphi}=\left\{ {\begin{array}{*{20}c}G ,\, if \, \left[ f_{0},f_{k-1} \right]>0,\\
-G ,\, if \, \left[ f_{0},f_{k-1} \right]<0.\\
\end{array} } \right. 
\end{equation}
\item If all poles of $Q \left( z \right)$ are real numbers and if  $\varphi^{1},\, \varphi^{2},\, \mathellipsis ,\,\varphi^{m}$ is the canonical set of all pole functions of $Q(z)$, i.e., root functions of $\hat{Q}(z)$, then the representing operator $A$ of $Q\left( z \right)$ in representation (\ref{eq312}) and the fundamental symmetry of the Pontryagin space $\mathcal{K}$ are given by the block matrices:
\begin{equation}
\label{eq321.2}
A=\left( {\begin{array}{*{20}c}
A_{\varphi^{1}} & 0 & ... & 0\\
0 & A_{\varphi^{2}} & ... & 0\\
\vdots & \vdots & \vdots & \vdots\\
0 & 0 & ... & A_{\varphi^{m}}\\
\end{array} } \right), \,
J=\left( {\begin{array}{*{20}c}
J_{\varphi^{1}} & 0 & ... & 0\\
0 & J_{\varphi^{2}} & ... & 0\\
\vdots & \vdots & \vdots & \vdots\\
0 & 0 & ... & J_{\varphi^{m}}\\
\end{array} } \right),
\end{equation}
respectively. The operator $A$ uniquely represents $\hat{L}(z)$ up to unitary isomorphism.
\end{enumerate}
\end{theorem}

\textbf{Proof}. 

(i) We will utilize results from \cite[p.76]{B1}. For this purpose we substitute $z_{1} $ with $\alpha$ and $\tilde{\varphi} (z)$ with $\varphi (z)$, meanwhile, the function $Q$, the representing operator $A$ of $Q$, and the Jordan chain (\ref{eq320}) of $A$ at $\alpha$ remain denoted here as in \cite{B1}.

According to \cite[p.76]{B1}, there exist a Jordan chain (\ref{eq320}) of the representing operator $A$ of $\hat{L}$ at $\alpha$ such that the limits (\ref{eq316}) and (\ref{eq318}) are satisfied.

The limit (\ref{eq318}) implies that the chain (\ref{eq320}) cannot be prolonged, i.e., it has a length $k$ equal to the order of the corresponding canonical root function $\varphi \left( z \right)$. 
Indeed, if the chain $f_{0},\, f_{1},\, \mathellipsis ,\, f_{k-1}$ of $A$ at $\alpha $ could be prolonged by a vector $f_{k}$, then $f_{k-1}=\left( A-\alpha \right)f_{k}$. It would follow 
\[
\left[ f_{0},f_{k-1} \right]=\left[ f_{0},\left( A-\alpha \right)f_{k}\right]=\left[ \left( A-\alpha \right)f_{0},f_{k} \right]=0,
\]
which contradicts (\ref{eq318}). 

Therefore, for every canonical root functions $\varphi$ of $\hat{Q}$ that satisfies limits (\ref{eq316}) and (\ref{eq318}), we can assign the subspace $S_{\varphi}$ of $ \mathcal{K}$ defined by (\ref{eq321}). The limit (\ref{eq318}) guarantees that the scalar product $\left[ .,.\right] $ does not degenerate on $S_{\varphi}$. It is well known, and easy to see, that Jordan vectors are linearly independent. This means that $S_{\varphi} $ is a Pontryagin space with $\dim S_{\varphi}=k< \infty $. Obviously, $S_{\varphi}$ is an invariant subspace of $A$. 

Recall, every operator in $\mathbb{C}^{k}$ has a matrix representation. We can select $f_{k-1} \in \mathbb{C}^{k}$ as 
\[
f_{k-1}=\left( 0,\, \mathellipsis ,0,\, 1 \right)^{T}.
\]
Since it holds 
\[
f_{k-2}=\left( A-\alpha \right)f_{k-1}=\left( 0,\, \mathellipsis ,1,\, 0 \right)^{T},\, \mathellipsis ,\, f_{0}=\left( 1,\, \mathellipsis ,0,\, 0 \right)^{T},
\]
the operator $A_{\varphi}=A_{\vert S_{\varphi}}:S_{\varphi}\to S_{\varphi}$ in matrix form must be the Jordan block (\ref{eq321.1}). Then we will define the indefinite scalar product by $\left[ x,y \right]_{S_{\varphi}}:=\left( J_{\varphi}x,y \right),\,  x,y\in S_{\varphi} \subseteq \mathbb{C}^{k}$ so that the Jordan block $A_{\varphi}$ will be self-adjoint with respect to that scalar product. That is equivalent to finding the fundamental symmetry $J_{\varphi}$ that satisfies $J_{\varphi}A_{\varphi}=A_{\varphi}^{T}J_{\varphi}$. Since it holds $GA_{\varphi}=A_{\varphi}^{T}G$, we can select either $J_{\varphi}=G$ or $J_{\varphi}=-G$. Because the scalar product $(J_{\varphi}.,.)$ must coincide with $\left[ x,y \right]_{\varphi}$ in $S_{\varphi}$, we will use the limit (\ref{eq318}) to determine $J_{\varphi}$. 

If we take $J_{\varphi}=G$, then we have 
\[
\left[ f_{0},f_{k-1} \right]=\left( \left( {\begin{array}{*{20}c}
0 & \cdots & 1\\
\vdots & \ddots & \vdots \\
1 & \cdots & 0\\
\end{array} } \right)\left( {\begin{array}{*{20}c}
1\\
\vdots \\
0\\
\end{array} } \right),\left( {\begin{array}{*{20}c}
0\\
\vdots \\
1\\
\end{array} } \right) \right)=\left( 1,1 \right)=1>0.
\]
If we select $J_{\varphi}=-G$, we have
\[
\left[ f_{0},f_{k-1} \right]=\left( \left( {\begin{array}{*{20}c}
0 & \cdots & -1\\
\vdots & \ddots & \vdots \\
-1 & \cdots & 0\\
\end{array} } \right)\left( {\begin{array}{*{20}c}
1\\
\vdots \\
0\\
\end{array} } \right),\left( {\begin{array}{*{20}c}
0\\
\vdots \\
1\\
\end{array} } \right) \right)=\left( -1,1 \right)=-1<0.
\]
Therefore, the sign of the limit (\ref{eq318}) uniquely determines the fundamental symmetry $J_{\varphi}$ of the Pontryagin space $S_{\varphi}$ as given by (\ref{eq321.15}). 

Because the scalar product $\left[.,.\right]_{\varphi}:=\left( J_{\varphi}.,.\right)$ does not degenerate on $S_{\varphi}$ and, according to limits (\ref{eq316}) and (\ref{eq318}), it coincides with scalar product of $\mathcal{K}$ on $S_{\varphi} $, the space $\left( S_{\varphi}, \left[ .,.\right]_{\varphi}  \right) $ is indeed a Pontryagin subspace of $\mathcal{K}$. 

(ii) By repeating the same process for all root functions $\varphi^{i}\left( z 
\right),\, i=1, 2, ..., m $, we will obtain the operator $A$ and the fundamental symmetry $J$ in the Jordan canonical form, as the matrix (\ref{eq321.2}) of the corresponding Jordan blocks on the diagonal. The matrices $A$ and $J$ are determined uniquely, up to order of the blocks on the diagonals of $A$ and $J$. 

It is easy to verify that $JA=A^{T}J$. This proves that operator $A$ is indeed a self-adjoint operator in the minimal Pontryagin space 
\begin{equation}
\label{eq321.3}
\mathcal{K}=S_{\varphi_{1}}[+] ... [+]S_{\varphi_{m}}.
\end{equation}

It remains to prove that the operator $A$ and the Krein space $\left( \mathcal{K}, \left[ .,. \right] \right)$ are unitarily isomorphic to any operator $\tilde{A}$ and a Pontryagin space $\tilde{\mathcal{K}}$  representing $Q(z)$ in a representation of the type (\ref{eq312}). 

As we explained in Section \ref{s2}.3, the poles, including orders, of a generalized Nevanlinna function $Q\in N_{\kappa}(\mathcal{H})$ represented by the minimal operator representation (\ref{eq12}), i.e., (\ref{eq312}), and poles of the resolvent $\left(A-z \right)^{-1} $ coincide. 

If in addition to representation (\ref{eq312}) we have a representation 
\begin{equation}
\label{eq321.4}
Q\left( z \right)=\tilde{S}+\tilde{\Gamma}^{+}\left(\tilde{A}-z \right)^{-1}\tilde{\Gamma} \in N_{\kappa 
}^{n\times n},\,  {\tilde{S}=\tilde{S}}^{\ast }\in \mathcal{L} \left( \mathbb{C}^{n} \right),
\end{equation}
then for the same canonical root function $\varphi \left( z \right)$ of $\hat{Q}(z)$ we have a Jordan chain $\tilde{f}_{0},\, \tilde{f}_{1},\, \mathellipsis $, $\tilde{f}_{k-1}$ of the operator $\tilde{A}$ at $\alpha $ that satisfies limits (\ref{eq316}) and (\ref{eq318}). Since the functions $Q\left( z \right)$ and $\varphi \left( z \right)$ are the same as before, the limits will remain the same. This further means that the order $k$ of the chain is the same as before, and that it holds $\left[ f_{0},f_{k-1} \right]=\left[ \tilde{f}_{0},\tilde{f}_{k-1} \right] $. Therefore, the mapping $U$ introduced by $U: \tilde{f}_{i} \rightarrow f_{i},\, i=1, 2, ..., k-1 $, is a unitary isomorphism between Pontryagin spaces $S_{\varphi}$ and  $ \tilde{S}_{\varphi}:=l.s. \left\{\tilde{ f}_{0},\, \tilde{f}_{1},\, \mathellipsis ,\, \tilde{f}_{k-1} \right\}$. Like $A$, the operator $\tilde{A}$ satisfies relations of the type $\tilde{f}_{i-1}=\left(\tilde{A}-\alpha \right)\tilde{f}_{i}, \, i=1, 2, ..., k-1 $. This means that the mapping $U$ is also an isomorphism between operators $A_{\vert S_{\varphi}} $ and $\tilde{A}_{\vert \tilde{S}_{\varphi}} $. 

By the same reasoning, we conclude that for every canonical root function $\varphi_{i}(z)$ spaces $\tilde{S}_{\varphi}$ and  $ S_{\tilde{\varphi}_{i}}$, as well as operators $A_{S_{\varphi_{i}}}$ and  $\tilde{A}_{S_{\tilde{\varphi}_{i}}}$ are unitarily isomorphic. 

This proves that the operator $A$ given by (\ref{eq321.2}) and the space $\mathcal{K} $ given by (\ref{eq321.3}) uniquely represent $Q(z)$ up to unitary isomorphism. \hfill $\square$

\begin{corollary}\label{corollary310} Let $L\left( z \right)$, $l\geq 1$, be an invertible Hermitian 
$n\times n, \, n\geq 2$, matrix polynomial that satisfies the condition \eqref{eq38}. Let $\varphi (z)$ be a canonical root function of $L$ at the eugenvalue $\alpha \in \mathbb{R}$  (obtained e.g., by (\ref{eq28})). Let $A$ be the representing operator of $\hat{L}$ and let the Jordan chain be given by \eqref{eq320}. Then the limits \eqref{eq316} and \eqref{eq318}, and all other statements of Theorem \ref{theorem38} hold when $\hat{L}$ and $\varphi$ replace $Q$ and $\varphi$.
\end{corollary} 

Assume that we have an invertible function $Q\in N_{\kappa }^{n \times n} $ and a pole $\alpha \in \mathbb{R}$ of $Q$. It is usually very difficult to obtain the inverse $\hat{Q}$. To avoid that difficulty, we may use the function $Q$ and its pole cancellation function $\psi (z)$ at $\alpha$ instead of the $\hat{Q}$ and $\varphi (z)$ in limits \eqref{eq316} and \eqref{eq318}. According to Definition \ref{defition24} we have
\begin{equation}
\label{eq322}
\varphi (z) =Q(z)\psi (z). 
\end{equation} 
Substituting this into limits \eqref{eq316} and \eqref{eq318}, we obtain the following statement:

\begin{theorem}\label{theorem311} Let $Q\in N_{\kappa }^{n \times n} $ satisfy the conditions of Proposition \ref{proposition32.1}, and let $A$ be the self-adjoint operator in the Pontryagin space $\left( \mathcal{K},\left[ .,.\right] \right) $ in the minimal representation (\ref{eq37}) of $Q$. 
\begin{enumerate}[(i)]
\item Let $\psi \left( z \right)$ be a canonical pole cancellation function of $Q$ at a pole $\alpha \in \mathbb{R}$ of order $k$ and $\varphi (z) =Q(z)\psi (z)$. Then there exists a Jordan chain \eqref{eq320} of the operator $A$ at $\alpha$ such that: 
\begin{equation}
\label{eq324}
\lim\limits_{z\to \alpha}{\frac{\psi \left( z \right)}{{(z-\alpha )}^{j}}}=0 \,\left( j < k \right),
\end{equation}
\begin{equation}
\label{eq326}
\lim\limits_{z\to \alpha}{\frac{\left(-\psi \left( z \right),Q(z)\psi \left( z \right) \right)}{{(z-\alpha )}^{k}}=\left[ f_{0},f_{k-1} \right]}\ne 0.
\end{equation}
The space $S_{\varphi}$ defined by \eqref{eq321}
is a Pontryagin subspace of $K$, and the operator $A_{\varphi}:=A_{\vert S_{\varphi}}$, with Jordan block form 
\eqref{eq321.1}, is a symmetric operator with respect to the indefinite scalar product $\left[ .,. \right]_{S_{\varphi}}:=\left[ .,. \right]_{\vert S_{\varphi}}=\left( J_{\varphi}.,.\right)  $.
The fundamental symmetry of the scalar product $\left[ .,. \right]_{S_{\varphi}}$ is given by 
\eqref{eq321.15}.
\item If all poles of $Q \left( z \right)$ are real numbers and if  $\psi^{1},\, \psi^{2},\, \mathellipsis ,\,\psi^{m}$ is the canonical set of all pole cancelation functions of $Q(z)$, then the representing operator $A$ of $Q\left( z \right)$ in representation (\ref{eq312}) and the fundamental symmetry of the Pontryagin space $\mathcal{K}$ are given by the block matrices \eqref{eq321.2}, respectively. The operator $A$ uniquely represents $\hat{L}(z)$ up to unitary isomorphism.
\end{enumerate}
\end{theorem}
The following example demonstrates one potential applications of the above theory in the area of the Hermitian matrix polynomials. 

\begin{example}\label{example312} Given $L\left( z \right)=\left( {\begin{array}{*{20}c}
1 & z^{2}-z\\
z^{2}-z & 0\\
\end{array} } \right)$.

\begin{enumerate}[(i)]
\item Find the eigenvalues of $L\left( z \right)$ and the corresponding canonical Jordan chains using Proposition \ref{proposition26}.
\item Find the operator representation of $\hat{L}\left( z \right)=\left( {\begin{array}{*{20}c}
0 & \frac{-1}{z\left( z-1 \right)}\\
\frac{-1}{z\left( z-1 \right)} & \frac{1}{{z^{2}\left( z-1 \right)}^{2}}\\
\end{array} } \right)$ using Theorem \ref{theorem38}.
\end{enumerate}
\end{example}
According to (\ref{eq22}), we will first multiply $L\left( z \right)$ by elementary matrices from the left 
and right until we obtain the diagonal form $D\left( z \right)$. It follows
\[
\left( {\begin{array}{*{20}c}
1 & 0\\
-\left( z^{2}-z \right) & 1\\
\end{array} } \right)\left( {\begin{array}{*{20}c}
1 & z^{2}-z\\
z^{2}-z & 0\\
\end{array} } \right)\left( {\begin{array}{*{20}c}
1 & -\left( z^{2}-z \right)\\
0 & 1\\
\end{array} } \right)
\]
\[
=\left( {\begin{array}{*{20}c}
1 & 0\\
0 & -\left( z^{2}-z \right)^{2}\\
\end{array} } \right)=:D\left( z \right).
\]
Expressing this in the notation of (\ref{eq22}), we have ${S\left( z \right)}\, 
L\left( z \right){T\left( z \right)}=D\left( z \right)$. According to 
(\ref{eq214}), for both eigenvalues $\alpha_{1}=0$ and $\alpha_{2}=1$, we have 
\[
P_{2}\left( 0 \right)=P_{2}\left( 1 \right)=\left( {\begin{array}{*{20}c}
0\\
1\\
\end{array} } \right)=:P_{2}.
\]
For both eigenvalues, according to (\ref{eq210}), we have
\[
\left( {\begin{array}{*{20}c}
1 & z^{2}-z\\
z^{2}-z & 0\\
\end{array} } \right)\left( {\begin{array}{*{20}c}
1 & -\left( z^{2}-z \right)\\
0 & 1\\
\end{array} } \right)\left( {\begin{array}{*{20}c}
0\\
1\\
\end{array} } \right)
\]
\[
=\left( {\begin{array}{*{20}c}
1 & 0\\
z^{2}-z & 1\\
\end{array} } \right)\left( {\begin{array}{*{20}c}
1 & 0\\
0 & -\left( z^{2}-z \right)^{2}\\
\end{array} } \right)\left( {\begin{array}{*{20}c}
0\\
1\\
\end{array} } \right).
\]
According to Corollary \ref{corollary210} the root functions in both eigrnvaues are equal 
\[
\varphi \left( z \right)=T\left( z \right)P_{2}=\left( 
{\begin{array}{*{20}c}
-z\left( z-1 \right)\\
1\\
\end{array} } \right).
\]
Note that for $Q(z)=\hat{L}(z)$ notation (\ref{eq322}) and (\ref{eq220}) are in sync. We have: 
\[
\psi \left( z \right)=L\left( z \right)\varphi \left( z \right)=\left( 
{\begin{array}{*{20}c}
1 & z^{2}-z\\
z^{2}-z & 0\\
\end{array} } \right)\left( {\begin{array}{*{20}c}
-\left( z^{2}-z \right)\\
1\\
\end{array} } \right)=\left( {\begin{array}{*{20}c}
0\\
-\left( z^{2}-z \right)^{2}\\
\end{array} } \right).
\]
Let us now focus on the eigenvalue $\alpha_{1}=0$ and denote the root function $\varphi^{1}(z)$. According to Proposition \ref{proposition26}, we have:
\begin{equation}
\label{eq328}
\varphi^{1} \left( z \right)=T\left( z \right)P_{2}=\left( 
{\begin{array}{*{20}c}
-z\left( z-1 \right)\\
1\\
\end{array} } \right)=\left( {\begin{array}{*{20}c}
0\\
1\\
\end{array} } \right)+\left( {\begin{array}{*{20}c}
1\\
0\\
\end{array} } \right)z-\left( {\begin{array}{*{20}c}
1\\
0\\
\end{array} } \right)z^{2}.
\end{equation}
Because $\alpha_{1}=0$ is a zero of the second order of the equation 
$d_{2}\left( z \right)=-z^{2}\left( z-1 \right)^{2}=0$, the maximal Jordan chain will have two Jordan vectors 
$\varphi^{1}_{0}$ and $\varphi^{1}_{1}$ at $\alpha_{1}=0$. According to (\ref{eq27})
and (\ref{eq328}) we conclude: 
\[
\varphi^{1}_{0}=\left( {\begin{array}{*{20}c}
0\\
1\\
\end{array} } \right)\, \wedge \varphi^{1}_{1}=\left( {\begin{array}{*{20}c}
1\\
0\\
\end{array} } \right)\wedge \tilde{\varphi}^{1}\left( z \right)=\left( 
{\begin{array}{*{20}c}
-1\\
0\\
\end{array} } \right).
\]
To confirm this conclusion, one can easily verify that vectors $\varphi^{1}_{0}$ and $\varphi^{1}_{1}$ satisfy the first two equations of system (\ref{eq120}) at $\alpha_{1}=0$, and that $\varphi^{1}_{0}$, $\varphi^{1}_{1}$, and $\tilde{\varphi}^{1}(\alpha_{1})$ do not satisfy the third equation of (\ref{eq120}). In other words, $\tilde{\varphi}^{1}\left( \alpha_{1} \right)=\left( {\begin{array}{*{20}c}
-1\\
0\\
\end{array} } \right)$ is not the third Jordan vector.

Let us now find the Jordan chain at the eigenvalue $\alpha_{2}=1$ and denote it by $\varphi^{2}(z)$. In this 
case, the expression (\ref{eq28}), $\varphi^{2} \left( z \right)=T\left( z \right)P_{2}$, is developed by exponents of $z-1$:
\[
\varphi^{2} \left( z \right)=\left( {\begin{array}{*{20}c}
-\left( z^{2}-z \right)\\
1\\
\end{array} } \right)=\left( {\begin{array}{*{20}c}
0\\
1\\
\end{array} } \right)+\left( {\begin{array}{*{20}c}
-1\\
0\\
\end{array} } \right)\left( z-1 \right)+\left( {\begin{array}{*{20}c}
-1\\
0\\
\end{array} } \right)\left( z-1 \right)^{2}.
\]
As before, we conclude $\varphi^{2}_{0}=\left( {\begin{array}{*{20}c}
0\\
1\\
\end{array} } \right)\, \wedge \varphi^{2}_{1}=\left( {\begin{array}{*{20}c}
-1\\
0\\
\end{array} } \right)\wedge \tilde{\varphi}^{2}\left( z \right)=\left( 
{\begin{array}{*{20}c}
-1\\
0\\
\end{array} } \right)$, 
where $\tilde{\varphi}^{2}\left( z \right)=\left( {\begin{array}{*{20}c}
-1\\
0\\
\end{array} } \right)$ is not a third Jordan vector of this chain.

Let us now calculate limit (\ref{eq326}) for both eigenvalues having in mind notation $Q(z)=\hat{L}(z)$.
\[
\frac{\left( \psi^{1}\left( z \right),\varphi^{1}\left( z \right) 
\right)}{\left( z-\alpha_{1} \right)^{k_{1}}}=\frac{\left( \left( 
{\begin{array}{*{20}c}
0\\
-\left( z^{2}-z \right)^{2}\\
\end{array} } \right),\left( {\begin{array}{*{20}c}
-z\left( z-1 \right)\\
1\\
\end{array} } \right) \right)}{z^{2}}=\frac{-z^{2}\left( z-1 
\right)^{2}}{z^{2}}\rightarrow-1\, (z\rightarrow 0).
\]
\[
\frac{\left( \psi^{2}\left( z \right),\varphi^{2}\left( z \right) 
\right)}{\left( z-\alpha_{2} \right)^{k_{2}}}=\frac{\left( \left( 
{\begin{array}{*{20}c}
0\\
-\left( z^{2}-z \right)^{2}\\
\end{array} } \right),\left( {\begin{array}{*{20}c}
-z\left( z-1 \right)\\
1\\
\end{array} } \right) \right)}{\left( z-1 \right)^{2}}=\frac{-z^{2}\left( 
z-1 \right)^{2}}{\left( z-1 \right)^{2}}\rightarrow-1\, (z\rightarrow 1).
\]
Let us denote by $A$ the representing operator of $\hat{L}\left( z 
\right)=-{L\left( z \right)}^{-1}$ and by $\lbrace f_{0}\left( 0 \right), 
f_{1}\left( 0 \right)\rbrace$ and $\lbrace f_{0}\left( 1 \right), f_{1}\left( 1 \right)\rbrace$ 
chains of $A$ at $\alpha_{1}=0$ and $\alpha_{2}=1$, respectively. 
According to (\ref{eq318}), we have:
\[
\left[ f_{0}\left( 0 \right),\, f_{1}\left( 0 \right) \right]=\left[ f_{0}\left( 1 
\right),\, f_{1}\left( 1 \right) \right]=-1.
\]
According to Theorem \ref{theorem311} (or Theorem \ref{theorem38}), the matrix $A$ and the fundamental symmetry 
$J$ will be:
\[
A=\left( {\begin{array}{*{20}c}
0 & 1 & 0 & 0\\
0 & 0 & 0 & 0\\
0 & 0 & 1 & 1\\
0 & 0 & 0 & 1\\
\end{array} } \right)\, \wedge \, J=\left( {\begin{array}{*{20}c}
0 & -1 & 0 & 0\\
-1 & 0 & 0 & 0\\
0 & 0 & 0 & -1\\
0 & 0 & -1 & 0\\
\end{array} } \right).
\]
It is easy to verify that $A$ is a self-adjoint operator in the Pontryagin space defined by $\left[ .,. \right]:=\left( J.,. \right)$. 
\
Now that we have $\hat{L}(z)$, the representing operator $A$, and the fundamental symmetry $J$, we can find operators $\Gamma $, and $\Gamma^{+} $, i.e., we can find the operator representation (\ref{eq312}). We have: 
\[
\left( A-z \right)^{-1}=\left( {\begin{array}{*{20}c}
\frac{-1}{z} & \frac{-1}{z^{2}} & 0 & 0\\
0 & \frac{-1}{z}& 0 & 0\\
0 & 0 & \frac{-1}{\left( z-1 \right)} & \frac{-1}{\left( z-1 \right)^{2}}\\
0 & 0 & 0 & \frac{-1}{\left( z-1 \right)}\\
\end{array} } \right).
\]
We will find $\Gamma $ from the equation 
\begin{equation}
\label{eq330}
\left( {\begin{array}{*{20}c}
0 & \frac{-1}{z\left( z-1 \right)}\\
\frac{-1}{z\left( z-1 \right)} & \frac{1}{{z^{2}\left( z-1 \right)}^{2}}\\
\end{array} } \right)=\Gamma^{+}\left( {\begin{array}{*{20}c}
\frac{-1}{z} & \frac{-1}{z^{2}} & 0 & 0\\
0 & \frac{-1}{z}& 0 & 0\\
0 & 0 & \frac{-1}{\left( z-1 \right)} & \frac{-1}{\left( z-1 \right)^{2}}\\
0 & 0 & 0 & \frac{-1}{\left( z-1 \right)}\\
\end{array} } \right)\Gamma. 
\end{equation}
Solving this equation for $\Gamma $ involves solving a system of nonlinear equations with eight unknowns $\gamma_{ij},\, i=1,2,3,4;\, j=1,2 $,  which is tedious but not difficult to solve. The solution is
\[
\Gamma =\left( {\begin{array}{*{20}c}
1 & 1\\
0 & 1\\
-1 & -1\\
0 & 1\\
\end{array} } \right)\Rightarrow \Gamma^{+}=\Gamma^{\ast }J=\left( 
{\begin{array}{*{20}c}
0 & -1 & 0 & 1\\
-1 & -1 & -1 & 1\\
\end{array} } \right).
\]
It is now easy to verify that $\hat{L}\left( z \right)=\Gamma^{+}\left( A-z \right)^{-1}\Gamma $.\hfill $\square$ 
\\

In the previous example, we were able to find operator $\Gamma$ from \eqref{eq330} because, in the representation (\ref{eq312}) of $\hat{L}$ it was $S=0$. In the case of a matrix polynomial $L\in N_{\kappa}^{n\times n}$, it is possible to have $S=\lim\limits_{z\to \infty}{{\hat{L}\left( z \right)}} \neq 0 $. An example is $L(z) =\left( {\begin{array}{*{20}c}
z & 1\\
1 & 1\\
\end{array} } \right) $. In that case, we would first find the limit $S$, and than solve a corresponding equation for $\tilde{\hat{L}}(z)=\hat{L}(z)-S$.
\\

\textbf{Compliance with Ethical Standards Statements:}

\textbf{Conflict of Interest:} The author declares that there is no conflict of interest. 

\textbf{Funding:} No funding was received to assist with the preparation of this manuscript.

\textbf{Ethical Conduct:}  This research is original. The manuscript is being submitted only to your journal. All previous results used in this paper are appropriately cited. 

\textbf{Data Availability Statements:} 

This research is purely mathematical. The data supporting this paper consist solely of mathematical statements, which are openly available in the cited references. There are no other data involved besides mathematical statements.

\end{document}